\documentclass[a4paper,12pt]{amsart}
\usepackage{amsfonts}
\usepackage{amssymb}
\usepackage{ifthen}
\usepackage{amscd}
\usepackage{amsxtra}
\usepackage{graphicx}
\usepackage{color}
\nonstopmode \numberwithin{equation}{section}
\newtheorem{thm}{Theorem}[section]
\newtheorem{lem}{Lemma}[section]
\newtheorem{cor}{Corollary}[section]
\newtheorem{prop}{Proposition}[section]

\newtheorem{step}{Step}[section]

\newtheorem{cl}{Claim}[section]
\newtheorem{ca}{Case}
\newtheorem{sca}{Subcase}[section]
\newtheorem{scl}{Subclaim}[section]
\newtheorem{conj}{Conjecture}

\theoremstyle{definition}
\newtheorem{defn}{Definition}[section]

\newtheorem{op}[equation]{Open Problem}
\newtheorem{ques}[equation]{Question}
\newtheorem{rem}{Remark}[section]
\newtheorem{exam}[equation]{Example}

\newcounter {own}
\def\theown {\thesection       .\arabic{own}}

\newenvironment{pf}[1][]{%
 \vskip 3mm
 \noindent
 \ifthenelse{\equal{#1}{}}%
  {{\slshape Proof. }}%
  {{\slshape #1.} }%
 }%
{\qed\bigskip}

\newcounter{alphabet}
\newcounter{tmp}

\makeatletter
\newcommand{\Ref}[1]{\@ifundefined{r@#1}{}{\setcounter{tmp}{\ref{#1}}\Alph{tmp}}}
\makeatother

\newenvironment{Lem}[1][]{\refstepcounter{alphabet}%
\bigskip%
\noindent%
{\bf Lemma \Alph{alphabet}}%
{\bf .} \itshape}{\vskip 8pt}

\newcommand{\Aut}{{\operatorname{Aut}}}

\def\be{\begin{equation}}
\def\ee{\end{equation}}

\newcommand{\bee}{\begin{enumerate}}
\newcommand{\eee}{\end{enumerate}}

\newcommand{\blem}{\begin{lem}}
\newcommand{\elem}{\end{lem}}
\newcommand{\bthm}{\begin{thm}}
\newcommand{\ethm}{\end{thm}}
\newcommand{\bcor}{\begin{cor}}
\newcommand{\ecor}{\end{cor}}
\newcommand{\beg}{\begin{exam}}
\newcommand{\eeg}{\end{exam}}
\newcommand{\begs}{\begin{examples}}
\newcommand{\eegs}{\end{examples}}
\newcommand{\bdefe}{\begin{defn}}
\newcommand{\edefe}{\end{defn}}
\newcommand{\bprob}{\begin{pro}}
\newcommand{\eprob}{\end{pro}}
\newcommand{\bques}{\begin{ques}}
\newcommand{\eques}{\end{ques}}
\newcommand{\bei}{\begin{itemize}}
\newcommand{\eei}{\end{itemize}}
\newcommand{\bcon}{\begin{conj}}
\newcommand{\econ}{\end{conj}}
\newcommand{\bop}{\begin{op}}
\newcommand{\eop}{\end{op}}

\newcommand{\bca}{\begin{ca}}
\newcommand{\eca}{\end{ca}}
\newcommand{\bsca}{\begin{sca}}
\newcommand{\esca}{\end{sca}}

\newcommand{\bcl}{\begin{cl}}
\newcommand{\ecl}{\end{cl}}

\newcommand{\bst}{\begin{step}}
\newcommand{\est}{\end{step}}

\newcommand{\bscl}{\begin{scl}}
\newcommand{\escl}{\end{scl}}

\newcommand{\bcons}{\begin{conjs}}
\newcommand{\econs}{\end{conjs}}
\newcommand{\bprop}{\begin{propo}}
\newcommand{\eprop}{\end{propo}}
\newcommand{\br}{\begin{rem}}
\newcommand{\er}{\end{rem}}
\newcommand{\brs}{\begin{rems}}
\newcommand{\ers}{\end{rems}}
\newcommand{\bo}{\begin{obser}}
\newcommand{\eo}{\end{obser}}
\newcommand{\bos}{\begin{obsers}}
\newcommand{\eos}{\end{obsers}}
\newcommand{\bpf}{\begin{pf}}
\newcommand{\epf}{\end{pf}}
\newcommand{\ba}{\begin{array}}
\newcommand{\ea}{\end{array}}
\newcommand{\beq}{\begin{eqnarray}}
\newcommand{\beqq}{\begin{eqnarray*}}
\newcommand{\eeq}{\end{eqnarray}}
\newcommand{\eeqq}{\end{eqnarray*}}

\newcommand{\ds}{\displaystyle}

\newcounter{minutes}
\divide\time by 60
\newcounter{hours}
\multiply\time by 60 \addtocounter{minutes}{-\time}

\begin{document}
\bibliographystyle{amsplain}
\title [] { Modulus of continuity and Heinz-Schwarz type inequalities of solutions to
 biharmonic equations}

\def\thefootnote{}
\footnotetext{ \texttt{\tiny File:~\jobname .tex,
          printed: \number\day-\number\month-\number\year,
          \thehours.\ifnum\theminutes<10{0}\fi\theminutes}
} \makeatletter\def\thefootnote{\@arabic\c@footnote}\makeatother


\author{Shaolin Chen${}^{~\mathbf{*}}$ }
\address{Sh. Chen, College of Mathematics and
Statistics, Hengyang Normal University, Hengyang, Hunan 421008,
People's Republic of China.} \email{mathechen@126.com}


\subjclass[2000]{Primary: 31A30, 31A05}
\keywords{Inhomogeneous  biharmonic Dirichlet problem,  modulus of
continuity, the Heinz-Schwarz type inequality.  \\
${}^{\mathbf{*}}$ Corresponding author (E-mail address:~
mathechen@126.com)}

\begin{abstract} 
For positive integers $n\geq2$ and $m\geq1$, consider a function $f$
satisfying the following: $(1)$ the inhomogeneous  biharmonic
equation $\Delta(\Delta f)=g$ ($g\in
\mathcal{C}(\overline{\mathbb{B}^{n}},\mathbb{R}^{m})$) in
$\mathbb{B}^{n}$, (2) the boundary conditions $f=\varphi_{1}$
$(\varphi_{1}\in \mathcal{C}(\mathbb{S}^{n-1},\mathbb{R}^{m}))$
 on $\mathbb{S}^{n-1}$ and $\partial
 f/\partial\mathbf{n}=\varphi_{2}$ ( $\varphi_{2}\in
\mathcal{C}(\mathbb{S}^{n-1},\mathbb{R}^{m})$)
  on $\mathbb{S}^{n-1}$, where $\partial /\partial\mathbf{n}$ stands for
the inward normal derivative, $\mathbb{B}^{n}$ is the unit ball in
$\mathbb{R}^{n}$ and $\mathbb{S}^{n-1}$ is the unit sphere of
$\mathbb{B}^{n}$. The main aim of this paper is to discuss the
Heinz-Schwarz type inequalities  and the modulus of continuity of
the solutions to the above inhomogeneous biharmonic Dirichlet
problem. 
\end{abstract}



\maketitle \pagestyle{myheadings} \markboth{Sh. Chen} {Modulus of
continuity and  Heinz-Schwarz type inequalities}

\section{ Introduction and  main results}\label{csw-sec1}
\subsection{Notations}
For a positive integer $n\geq2$, let $\mathbb{R}^{n}$ and
$\mathbb{R}=\mathbb{R}^{1}$ be
  the usual real vector space of dimension $n$ and  the set of real numbers, respectively.
Let $\mathbb{B}^{n}(x_{0},r)=\{x\in\mathbb{R}^{n}:~|x-x_{0}|<r\}$,
$\overline{\mathbb{B}^{n}}(x_{0},r)=\{x\in\mathbb{R}^{n}:~|x-x_{0}|\leq
r\}$ and
$\mathbb{S}^{n-1}(x_{0},r)=\partial\mathbb{B}^{n}(x_{0},r)$, where
$r>0$. We write $\mathbb{B}^{n}:=\mathbb{B}^{n}(0,1)$ and
$\mathbb{S}^{n-1}:=\mathbb{S}^{n-1}(0,1)$.
 Set $\mathbb{D}=\mathbb{B}^2$, the
open unit disk in the complex plane $\mathbb{C}\cong
\mathbb{R}^{2}$. For $m\in\mathbb{N}:=\{1,2,\ldots\}$ and
$k\in\mathbb{N}_0=\mathbb{N}\cup\{0\}$, we denote by
$\mathcal{C}^{k}(\Omega_{1},\Omega_{2})$ the set of all $k$-times
continuously
 differentiable functions from $\Omega_{1}$ into
$\Omega_{2}$, where $\Omega_{1}$ and $\Omega_{2}$ are subsets of
$\mathbb{R}^{n}$ and $\mathbb{R}^{m}$, respectively. In particular,
let
$\mathcal{C}(\Omega_{1},\Omega_{2}):=\mathcal{C}^{0}(\Omega_{1},\Omega_{2})$,
the set of all continuous functions of $\Omega_{1}$ into
$\Omega_{2}$.

\subsection{Inhomogeneous  biharmonic
equation} For $n\geq2$ and $m\geq1$, let  $\varphi_{1}\in
\mathcal{C}(\mathbb{S}^{n-1},\mathbb{R}^{m})$, $\varphi_{2}\in
\mathcal{C}(\mathbb{S}^{n-1},\mathbb{R}^{m})$ and $g\in
\mathcal{C}(\overline{\mathbb{B}^{n}},\mathbb{R}^{m})$. Of
particular interest to us is the following {\it inhomogeneous
biharmonic problem}:

\be\label{eq-c-1}
\begin{cases}
\displaystyle \Delta(\Delta f)=g&\mbox{ in }\, \mathbb{B}^{n},\\
\displaystyle  f=\varphi_{1} &\mbox{ on }\, \mathbb{S}^{n-1},\\
\displaystyle \frac{\partial
f}{\partial\mathbf{n}}=\varphi_{2}&\mbox{ on }\, \mathbb{S}^{n-1},
\end{cases}
\ee where $\Delta:=\sum_{j=1}^{n}\frac{\partial^{2}}{\partial
x_{j}^{2}} $ is the Laplace operator,
$$\frac{\partial f}{\partial\mathbf{n}}=\left(\frac{\partial
f_{1}}{\partial\mathbf{n}},\ldots,\frac{\partial
f_{m}}{\partial\mathbf{n}}\right)$$ and $\partial
f_{k}/\partial\mathbf{n}$ denotes the differentiation in the inward
normal direction for $k\in\{1,\ldots,m\}$.
 Here the boundary
conditions in  (\ref{eq-c-1}) are interpreted in the following
distributional sense. For some fixed $r\in(0,1)$, let
$F(x)=f(rx),~x\in\mathbb{B}^{n}$. Then

\be\label{B-1}
\begin{cases}
\displaystyle \Delta(\Delta F(x))=r^{2}\Delta(\Delta f(rx)),& x\in\mathbb{B}^{n},\\
\displaystyle  F=f_{r} &\mbox{ on }\, \mathbb{S}^{n-1},\\
\displaystyle
 \frac{\partial
F}{\partial\mathbf{n}}=r\frac{\partial
f_{r}}{\partial\mathbf{n}}&\mbox{ on }\, \mathbb{S}^{n-1},
\end{cases}
\ee and
 $f_{r}\rightarrow\varphi_{1}$ as $r\rightarrow1^{-}$, and
$r\big(\partial f_{r}/\partial\mathbf{n}\big)\rightarrow\varphi_{2}$
as $r\rightarrow1^{-}$, where
$$f_{r}(\zeta):=f(r\zeta),~\zeta\in\mathbb{S}^{n-1}.$$

In particular, if $g\equiv0$, then the solutions to (\ref{eq-c-1})
are {\it biharmonic mappings} (see
\cite{CLW-2018,C-Z,GMR-2018,S-2003}).




The  inhomogeneous biharmonic equations arise in areas of continuum
mechanics, including linear elasticity theory and the solution of
Stokes flows (cf. \cite{Ha,Kh-1996,Se,We}). 
 This article continues the study of the previous
work of Li et al. \cite{LP},  Kalaj \cite{K5} and the monograph of
Gazzola et al. \cite{GGS}. In order to state our main results, we
introduce some necessary terminologies.

For $x\in\mathbb{B}^{n}$ and $\zeta\in\mathbb{S}^{n-1}$, let
\beqq\label{eq-c-3}H_{n}(x,\zeta)=\frac{1}{2}\frac{(1-|x|^{2})^{2}}{|x-\zeta|^{n}}\eeqq
and
\beqq\label{eq-c-4}K_{n}(x,\zeta)=\frac{1}{4}\frac{(1-|x|^{2})^{2}}{|x-\zeta|^{n+2}}\big(n(1-|x|^{2})-(n-4)|x-\zeta|^{2}\big).
\eeqq Here $K_{n}$ is called a {\it  biharmonic  Poisson kernel}
(see \cite[p.157]{GGS}).

For  $x,y\in\mathbb{R}^{n}\backslash\{0\}$, we define
$x^{\ast}=x/|x|$, $y^{\ast}=y/|y|$,
$$[x,y]:=\left|y|x|-x^{\ast}\right|=\left|x|y|-y^{\ast}\right|~
\mbox{and} ~x\otimes y:=(1-|x|^{2})(1-|y|^{2}).$$ Also, for
$x,y\in\mathbb{B}^{n}$ with $x\neq y$, we use $G_{2,n}(x,y)$ to
denote the {\it biharmonic Green function}:

 \be\label{eq-c-5} G_{2,n}(x,y)=
\begin{cases}
\displaystyle
c_{n}\left(|x-y|^{4-n}-[x,y]^{4-n}-\frac{(n-4)}{2}(x\otimes
y)[x,y]^{2-n}\right),
 &n\neq2,4,\\
 \displaystyle\\
 \displaystyle c_{n}\left(|x-y|^{4-n}\log\frac{|x-y|^{2}}{[x,y]^{2}} +
 (x\otimes y)[x,y]^{2-n}\right),   &n=2,4,
\end{cases}\ee
where $A_{n-1}=2\pi^{\frac{n}{2}}/\Gamma\big(\frac{n}{2}\big)$ is
the $(n-1)$-dimensional surface area of $\mathbb{S}^{n-1}$,
$c_{n}=1/\big(2(4-n)(2-n)A_{n-1}\big)$ for $n\neq2,4$, and
$c_{n}=1/\big(8(-1)^{n/2+1}A_{n-1}\big)$ for $n=2,4$.  We refer
readers to the Chapter 4 in \cite{GGS} for general properties of the
Green functions.

\subsection{Main results} For $n=2$, Li and Ponnusamy (\cite[Theorem 1.1]{LP}) established a
representation formula and the uniqueness of the solutions to
(\ref{eq-c-1}). In fact, for $n\geq2$, the solutions to
(\ref{eq-c-1}) with smooth boundary conditions   has alreadly been
observed in \cite[p.138]{GGS}. For the sake of completeness, we
recall
 the representation formula and the
uniqueness of the solutions to (\ref{eq-c-1}) for some slightly
weaker boundary value conditions.

\begin{prop}\label{thm-1}
For positive integers $n\geq2$ and $m\geq1$, suppose that
$\varphi_{1}\in \mathcal{C}(\mathbb{S}^{n-1},\mathbb{R}^{m})$,
$\varphi_{2}\in \mathcal{C}(\mathbb{S}^{n-1},\mathbb{R}^{m})$ and
$g\in \mathcal{C}(\overline{\mathbb{B}^{n}},\mathbb{R}^{m})$. Let

$$K[\varphi_{1}](x)=\int_{\mathbb{S}^{n-1}}K_{n}(x,\zeta)\varphi_{1}(\zeta)d\sigma(\zeta),~H[\varphi_{2}](x)=\int_{\mathbb{S}^{n-1}}H_{n}(x,\zeta)\varphi_{2}(\zeta)d\sigma(\zeta)$$
and $$G[g](x)=\int_{\mathbb{B}^{n}}G_{2,n}(x,y)g(y)dV(y),$$ where
$d\sigma$ denotes the normalized Lebesgue surface measure on
$\mathbb{S}^{n-1}$ and $dV$ is the Lebesgue volume measure on
$\mathbb{B}^{n}$. If $f$ is a solution to (\ref{eq-c-1}), then

$$f(x)=K[\varphi_{1}](x)+H[\varphi_{2}](x)+G[g](x),~x\in\mathbb{B}^{n}.$$
\end{prop}

Heinz in his classical paper \cite{He} showed that the following
result which is called the Heinz-Schwarz type inequality  of
harmonic mappings: If $f$ is a harmonic mapping of $\mathbb{D}$ into
$\mathbb{D}$ with $f(0)=0$, then
$$|f(z)|\leq\frac{4}{\pi}\arctan|z|, ~z\in\mathbb{D}.$$ Later, Hethcote
 \cite{Het} removed
the assumption $f(0)=0$ and proved the following inequality

\beqq\label{eq-pav1}\left|f(z)-\frac{1-|z|^{2}}{1+|z|^{2}}f(0)\right|\leq\frac{4}{\pi}\arctan
|z|,~z\in\mathbb{D},\eeqq where $f$ is a harmonic mapping from
$\mathbb{D}$ into itself (see also  \cite[Theorem 3.6.1]{Pav1}). For
$n\geq3$, the classical Heinz-Schwarz type inequality of harmonic
mappings in $\mathbb{B}^{n}$ infers that if $f$ is a harmonic
mapping of $\mathbb{B}^{n}$ into itself satisfying $f(0)=0,$ then
$$|f(x)|\leq U(|x|e_{n}),
$$
where  $e_{n}=(0,\ldots,0,1)$ and $U$ is a harmonic function of
$\mathbb{B}^{n}$ into $[-1,1]$ defined by
$$U(x)=P[\chi_{\mathbb{S}^{n-1}_{+}}-\chi_{\mathbb{S}^{n-1}_{-}}](x):=
\int_{\mathbb{S}^{n-1}}\frac{1-|x|^{2}}{|x-\zeta|^{n}}\left(\chi_{\mathbb{S}^{n-1}_{+}}(\zeta)-\chi_{\mathbb{S}^{n-1}_{-}}(\zeta)\right)d\sigma(\zeta).
$$
Here $\chi$ is the indicator function and
$\mathbb{S}^{n-1}_{+}=\{x\in\mathbb{S}^{n-1}:~x_{n}\geq0\}$,
$\mathbb{S}^{n-1}_{-}=\{x\in\mathbb{S}^{n-1}:~x_{n}\leq0\}$ (cf.
\cite{ABR}). In \cite{K5}, Kalaj showed  the following result for
harmonic mappings $f$ of $\mathbb{B}^{n}$ into itself:
\be\label{eq-K}
\left|f(x)-\frac{1-|x|^{2}}{(1+|x|^{2})^{\frac{n}{2}}}f(0)\right|\leq
U(|x|e_{n}),~x\in\mathbb{B}^{n}.\ee

By analogy with the inequality (\ref{eq-K}), we obtain the following
result.

\begin{thm}\label{thm-3} For positive integers $n\geq2$ and $m\geq1$,
 let $\varphi_{1}\in \mathcal{C}(\mathbb{S}^{n-1},\mathbb{R}^{m})$,
$\varphi_{2}\in \mathcal{C}(\mathbb{S}^{n-1},\mathbb{R}^{m})$ and
$g\in \mathcal{C}(\overline{\mathbb{B}^{n}},\mathbb{R}^{m})$. If $f$
is a solution to (\ref{eq-c-1}), then, for
$x\in\overline{\mathbb{B}^{n}}$,

\beq\label{eq-th-3} \big|f(x)&-&\delta_{1}(|x|)K[\varphi_{1}](0)-\delta_{2}(|x|)H[\varphi_{2}](0)\big|\\
\nonumber
&\leq&\frac{|n-4|}{4}\|\varphi_{1}\|_{\infty}(1-|x|^{2})U(|x|e_{n})
+\frac{n}{4}\|\varphi_{1}\|_{\infty}U^{\ast}(|x|e_{n})\\ \nonumber
&&+\frac{\|\varphi_{2}\|_{\infty}}{2}(1-|x|^{2})U(|x|e_{n})+\frac{\|g\|_{\infty}}{8n(n+2)}(1-|x|^{2})^{2},\eeq
where
$\|\varphi_{k}\|_{\infty}=\sup_{\zeta\in\mathbb{S}^{n-1}}|\varphi_{k}(\zeta)|~(k=1,2),$
$\|g\|_{\infty}=\sup_{x\in\mathbb{B}^{n}}|g(x)|,$

$$\delta_{1}(|x|)=\frac{n}{4}\frac{(1-|x|^{2})^{3}}{(1+|x|^{2})^{\frac{n+2}{2}}}-\frac{(n-4)}{4}\frac{(1-|x|^{2})^{2}}{(1+|x|^{2})^{\frac{n}{2}}},
~\delta_{2}(|x|)=\frac{(1-|x|^{2})^{2}}{(1+|x|^{2})^{\frac{n}{2}}}$$
and
$$U^{\ast}(x)=\int_{\mathbb{S}^{n-1}}\frac{(1-|x|^{2})^{3}}{|x-\zeta|^{n+2}}\left(\chi_{\mathbb{S}^{n-1}_{+}}(\zeta)-\chi_{\mathbb{S}^{n-1}_{-}}(\zeta)\right)d\sigma(\zeta).$$

Moreover, if  $g=(8n(2+n)M,0,\ldots,0)\in\mathbb{R}^{m}$ in
$\mathbb{B}^{n}$, and
$\varphi_{1}=\varphi_{2}=(0,\ldots,0)\in\mathbb{R}^{m}$ in
$\mathbb{S}^{n-1}$, then
$f(x)=\left(M(1-|x|^{2})^{2},0,\ldots,0\right)\in\mathbb{R}^{m}$
shows that the estimate of (\ref{eq-th-3}) is sharp in
$\overline{\mathbb{B}^{n}}$, where $M$ is a constant.

\end{thm}


We remark that Theorem \ref{thm-3} is somehow weakened by the fact
that the function $\delta_{1}(|x|)$ may change sign; this is in
contrast with what happens in the second order case (see
(\ref{eq-K})).

\begin{table}
\center{
\begin{tabular}{|c|c||c|c|c|c|}
\hline $n$ & $U^{\ast}(re_{n})$ & $U(re_{n})$
\\
\hline 2&
$\frac{2}{\pi}\left(\frac{2r(1-r^{2})}{1+r^{2}}+(1+r^{2})\arcsin\frac{2r}{1+r^{2}}\right)$
& $\frac{4 \arctan r}{\pi}$
\\
3&
$\frac{\left((1+4r^{2}+3r^{4})(1+r^{2})^{\frac{1}{2}}+r^{6}-3r^{4}+3r^{2}-1\right)}{3
r(1+r^{2})^{\frac{3}{2}}}$ &
$\frac{r^{2}-1+\sqrt{1+r^{2}}}{r\sqrt{1+r^{2}}}$
\\
4&
$\frac{2}{\pi}\left(\frac{2r(1-r^{2})}{(1+r^{2})^{2}}+\arcsin\frac{2r}{1+r^{2}}\right)$
&$\frac{2r(r^{2}-1)+2(1+r^{2})^{2}\arctan r}{\pi
r^{2}(1+r^{2})}$\\
\hline
\end{tabular}
}
\bigskip
\caption{Values of $U^{\ast}(re_{n})$ and $U(re_{n})$ for Theorem
\ref{thm-3}.}
\end{table}





A continuous increasing function $\omega:\, [0,\infty)\rightarrow
[0,\infty)$ with $\omega(0)=0$ is called a {\it majorant} if
$\omega(t)/t$ is non-increasing for $t>0$ (cf. \cite{Dy1,Dy2}).
Given a subset $\Omega$ of $\mathbb{R}^{n}$, a function
$u:~\Omega\rightarrow\mathbb{R}^{m}$ is said to belong to the {\it
Lipschitz space} $\mathcal{L}_{\omega}(\Omega,\mathbb{R}^{m})$ if
there is a positive constant $L$ such that
$$\sup_{x,y\in\Omega, x\neq y}\frac{|u(x)-u(y)|}{\omega(|x-y|)}\leq L.$$

Dyakonov \cite{Dy1,Dy2} characterized the analytic functions of
class $\mathcal{L}_{\omega}(\mathbb{D},\mathbb{C})$
 in terms of their modulus (see also \cite{Pav}).

 It is well-known that the condition $u\in
\mathcal{L}_{\omega}(\mathbb{S}^{1},\mathbb{C})$ is not enough to
guarantee that its harmonic extension $P[\psi]$ belongs to
$\mathcal{L}_{\omega}(\mathbb{D},\mathbb{C})$, where $\omega(t)=t$
and
$$P[\psi](z)=\frac{1}{2\pi}\int_{0}^{2\pi}\frac{1-|z|^{2}}{|e^{i\theta}-z|^{2}}u(e^{i\theta})d\theta,~z\in\mathbb{D}.$$ In fact,
$P[\psi]\in\mathcal{L}_{\omega}(\mathbb{D},\mathbb{C})$ is Lipschitz
continuous if and only if the {\it Hilbert transform} of
$d\psi(e^{i\theta})/d\theta$ belongs to $L^{\infty}(\mathbb{S}^{1})$
(see \cite{AKM} and \cite{Z}), where $\omega(t)=t$. In \cite{AKM},
Arsenovi\'c et al. established the following result for harmonic
mappings of $\mathbb{B}^n$ into $\mathbb{R}^{n}$: For a boundary
function which is Lipschitz continuous, if its harmonic extension is
{\it quasiregular}, then this extension is also Lipschitz
continuous. Recently, the relationship of the Lipschitz continuity
between the boundary functions and their harmonic extensions has
attracted much attention (see \cite{CHRW,CLW-2018,K2,K3,LP}).  Li
and Ponnusamy \cite{LP} discussed the Lipschitz characteristic of
solutions to the inhomogeneous biharmonic equation (\ref{eq-c-1})
for $n=2$.  The same problem in higher dimentional space is much
more complicated because of the lack of the techniques of complex
analysis. For $n\geq2$, we will investigate the Lipschitz continuity
(or the modulus of continuity) of the solutions to (\ref{eq-c-1}) as
follows.


\begin{thm}\label{thm-2}
Suppose that   $n\geq2$ and $m\geq1$ are  integers, and $\omega$ is
a majorant satisfying
$$\limsup_{t\rightarrow0^{+}}\frac{\omega(t)}{t}=c<\infty.$$ For $\varphi_{1}\in
\mathcal{L}_{\omega}(\mathbb{S}^{n-1},\mathbb{R}^{m})$,
$\varphi_{2}\in \mathcal{C}(\mathbb{S}^{n-1},\mathbb{R}^{m})$ and
$g\in \mathcal{C}(\overline{\mathbb{B}^{n}},\mathbb{R}^{m}),$ if $f$
satisfies (\ref{eq-c-1}), then
$f\in\mathcal{L}_{\omega}(\mathbb{B}^{n},\mathbb{R}^{m}).$
\end{thm}

The rest of this article is organized as follows. In section
\ref{csw-sec2}, some necessary notations and useful results will be
introduced. In section \ref{csw-sec5}, the Proposition \ref{thm-1}
and Theorem \ref{thm-3} will be proved. Theorem \ref{thm-2} will be
showed in section \ref{csw-sec4}.



\section{ Preliminaries}\label{csw-sec2}
\subsection{Gauss Hypergeometric Functions}\label{sbcsw-sec2.1}

For $a, b, c\in\mathbb{R}$ with $c\neq0, -1, -2, \ldots,$ the {\it
hypergeometric} function is defined by the power series
$${_{2}}F_{1}(a,b;c;t)=\sum_{k=0}^{\infty}\frac{(a)_{k}(b)_{k}}{(c)_{k}}\frac{t^{k}}{k!}
$$ with respect
to the variable $t\in(-1,1)$. Here $(a)_{0}=1$,
$(a)_{k}=a(a+1)\cdots(a+k-1)$ for $k=1, 2, \ldots$, and generally
$(a)_{k}=\Gamma(a+k)/\Gamma(a)$ is the {\it Pochhammer} symbol,
where $\Gamma$ is the {\it Gamma function} (cf. \cite{PBM}).

\subsection{M\"obius Transformations of $\mathbb{B}^{n}$}\label{sbcsw-sec2.2}
For any fixed $x\in\mathbb{B}^{n}$, the {\it M\"obius
transformation} in $\mathbb{B}^{n}$ is defined by \be\label{eq-ex1}
\phi_{x}(y)=\frac{|x-y|^{2}x-(1-|x|^{2})(y-x)}{[x,y]^{2}},~y\in\mathbb{B}^{n}.
\ee The set of isometries of the hyperbolic unit ball is a {\it
Kleinian subgroup} of all M\"obius transformations of the extended
spaces $\mathbb{R}^{n}\cup\{\infty\}$ onto itself. In the following,
we make use of the {\it automorphism group} $\Aut(\mathbb{B}^{n})$
consisting of all M\"obius transformations of the unit ball
$\mathbb{B}^{n}$ onto itself. We recall the following facts from
\cite{Bea}: For $x\in\mathbb{B}^{n}$ and
$\phi_{x}\in\Aut(\mathbb{B}^{n})$, we have $\phi_{x}(0)=x$,
$\phi_{x}(x)=0$, $\phi_{x}(\phi_{x}(y))=y \in\mathbb{B}^{n}$,

\be\label{II} |\phi_{x}(y)|=\frac{|x-y|}{[x,y]},
~1-|\phi_{x}(y)|^{2}=\frac{(1-|x|^{2})(1-|y|^{2})}{[x,y]^{2}}\ee and
\be\label{III} |J_{\phi_{x}}(y)|=\frac{(1-|x|^{2})^{n}}{[x,y]^{2n}}.
\ee

\subsection{Matrix notations}\label{sbcsw-sec2.3}

For an $m\times n$ matrix $A=(a_{ij})_{m\times n}$, the operator
norm of $A$ is defined by
$$|A|=\sup_{x\neq 0}\frac{|Ax|}{|x|}=\max\{|A\theta|:\,
\theta\in\mathbb{S}^{n-1}\},
$$
and the matrix function $l(A)$ is defined by
$$l(A)=\inf\{|A\theta|:~\theta\in\mathbb{S}^{n-1}\}.$$

For a domain $\Omega\subset\mathbb{R}^{n}$, let
$f=(f_{1},\ldots,f_{m}):~\Omega\rightarrow\mathbb{R}^{m}$ be a
function that has all partial derivatives at
$x=(x_{1},\ldots,x_{n})$ in $\Omega$. Then we denote the derivative
$D_{f}$ of $f$ by
$$D_{f}=\left(\begin{array}{cccc}
\ds D_{1}f_{1}\; \cdots\;
 D_{n}f_{1}\\[4mm]
\vdots\;\; \;\;\cdots\;\;\;\;\vdots \\[2mm]
 \ds
D_{1}f_{m}\; \cdots\;
 D_{n}f_{m}
\end{array}\right)=(\nabla f_{1}~\cdots,\nabla f_{m})^{T},
$$ where $D_{j}f_{i}(x)=\partial f_{i}(x)/\partial x_j$, $T$ is the
transpose and the gradients $\nabla f_{j}~(j=1,\ldots,m)$ are
understood as column vectors.

\subsection{Spherical coordinate transformation}\label{sbcsw-sec2.4}

Let
$Q=(\zeta_{1},\ldots,\zeta_{n}):~\mathbb{T}^{n-1}\rightarrow\mathbb{S}^{n-1}$
be the following spherical coordinate transformation
\begin{eqnarray*}
\zeta_{1}&=&\cos\theta_{1},\\
\zeta_{2}&=&\sin\theta_{1}\sin\theta_{2},\\
&\vdots&\\
\zeta_{n-1}&=&\sin\theta_{1}\sin\theta_{2}\cdots\sin\theta_{n-2}\cos\theta_{n-1},\\
\zeta_{n}&=&\sin\theta_{1}\sin\theta_{2}\cdots\sin\theta_{n-2}\sin\theta_{n-1},
\end{eqnarray*}
where
$\mathbb{T}^{n-1}=\underbrace{[0,\pi]\times\cdots\times[0,\pi]}\times[0,2\pi].$
We use
$$J_{Q}(\theta_{1},\ldots,\theta_{n-1}):=\sin^{n-2}\theta_{1}\cdots\sin\theta_{n-2}$$ to denote the {\it Jacobian} of $Q$.


\section{The Heinz-Schwarz type inequalities of solutions to inhomogeneous
 biharmonic Dirichlet problems}\label{csw-sec5}

\subsection*{The proof of Proposition \ref{thm-1}}
For some fixed $r\in(0,1)$, let $F(x)=f(rx),~x\in\mathbb{B}^{n}$.
  It follows from \cite[Formula 4.98]{GGS}
that
$$F(x)=K[\varphi_{1}^{\ast}](x)+H[\varphi_{2}^{\ast}](x)+G[g^{\ast}](x)$$ is the only solution to (\ref{B-1}), where  $g^{\ast}(x)=r^{2}g(rx),~x\in\mathbb{B}^{n}$, $\varphi_{1}^{\ast}(\zeta)=f_{r}(\zeta):=f(r\zeta)$ and
$\varphi_{2}^{\ast}(\zeta)=r\frac{\partial
f_{r}(\zeta)}{\partial\mathbf{n} }, ~\zeta\in\mathbb{S}^{n-1}$. By
letting $r\rightarrow1^{-}$, we get the desired result. \qed

The following result will be used in the proof of Theorem
\ref{thm-3}.

\begin{Lem}{\rm (\cite{K-2017}~$\mbox{or}$~\cite[2.5.16(43)]{PBM})}\label{pro-1}
For $\mu_{1}>1$ and $\mu_{2}>0$, we have
$$\int_{0}^{\pi}\frac{\sin^{\mu_{1}-1}t}{(1+r^{2}-2r\cos t)^{\mu_{2}}}dt=
\mathbf{B}\left(\frac{\mu_{1}}{2},\frac{1}{2}\right)
{_{2}}F_{1}\big(\mu_{2},\mu_{2}+\frac{1-\mu_{1}}{2};\frac{1+\mu_{1}}{2};r^{2}\big),~r\in[0,1),
$$
where $\mathbf{B}(.,.)$ denotes the beta function.
\end{Lem}

\subsection*{The proof of Theorem \ref{thm-3}} For $x\in\mathbb{B}^{n}$, let
$$\Phi(x)=K[\varphi_{1}](x)+H[\varphi_{2}](x).
$$ Then $\Phi$ is biharmonic in $\mathbb{B}^{n}$.

We first assume that $x=|x|e_{n}$ is on the ray $[0,e_{n}]$, where
$e_{n}=(0,\ldots,0,1)$. Then we have

\be\label{eq-xj-5}
\big|\Phi(|x|e_{n})-\delta_{1}(|x|)K[\varphi_{1}](0)-\delta_{2}(|x|)H[\varphi_{2}](0)
\big|\leq|K^{\ast}(|x|e_{n})|+|H^{\ast}(|x|e_{n})|,
 \ee
where
$$\delta_{1}(|x|)=\frac{n}{4}\frac{(1-|x|^{2})^{3}}{(1+|x|^{2})^{\frac{n+2}{2}}}-\frac{(n-4)}{4}\frac{(1-|x|^{2})^{2}}{(1+|x|^{2})^{\frac{n}{2}}},
~\delta_{2}(|x|)=\frac{(1-|x|^{2})^{2}}{(1+|x|^{2})^{\frac{n}{2}}},$$
$$K^{\ast}(|x|e_{n})=K[\varphi_{1}](|x|e_{n})-\delta_{1}(|x|)K[\varphi_{1}](0)$$
and
$$H^{\ast}(|x|e_{n})=H[\varphi_{2}](|x|e_{n})-\delta_{2}(|x|)H[\varphi_{2}](0).$$

By calculations, we obtain  \beq \label{eq-xj-6}\left|
K^{\ast}(|x|e_{n})\right|&\leq&\frac{n}{4}\|\varphi_{1}\|_{\infty}\int_{\mathbb{S}^{n-1}}
\left|\frac{(1-|x|^{2})^{3}}{||x|e_{n}-\zeta|^{n+2}}-\frac{(1-|x|^{2})^{3}}{(1+|x|^{2})^{\frac{n+2}{2}}}\right|d\sigma(\zeta)\\
\nonumber
&&+\frac{|n-4|}{4}\|\varphi_{1}\|_{\infty}\int_{\mathbb{S}^{n-1}}
\left|\frac{(1-|x|^{2})^{2}}{||x|e_{n}-\zeta|^{n}}-\frac{(1-|x|^{2})^{2}}{(1+|x|^{2})^{\frac{n}{2}}}\right|d\sigma(\zeta)\\
\nonumber
&=&\frac{n}{4}\|\varphi_{1}\|_{\infty}U^{\ast}(|x|e_{n})+\frac{|n-4|}{4}\|\varphi_{1}\|_{\infty}(1-|x|^{2})U(|x|e_{n})
\eeq and

 \beq\label{eq-xj-7} |H^{\ast}(|x|e_{n})|&\leq&
\frac{1}{2}\|\varphi_{2}\|_{\infty}\int_{\mathbb{S}^{n-1}}\bigg|\frac{(1-|x|^{2})^{2}}{|x-\zeta|^{n}}-\frac{(1-|x|^{2})^{2}}{(1+|x|^{2})^{\frac{n}{2}}}\bigg|d\sigma(\zeta)\\
\nonumber
&\leq&\frac{(1-|x|^{2})}{2}\|\varphi_{2}\|_{\infty}U(|x|e_{n}).\eeq

It follows from (\ref{eq-xj-5}), (\ref{eq-xj-6}) and (\ref{eq-xj-7})
that

\beq\label{eq-xj-8}
|\Phi^{\ast}(|x|e_{n})|&\leq&\frac{|n-4|}{4}\|\varphi_{1}\|_{\infty}(1-|x|^{2})U(|x|e_{n})
+\frac{n}{4}\|\varphi_{1}\|_{\infty}U^{\ast}(|x|e_{n})\\
\nonumber
&&+\|\varphi_{2}\|_{\infty}\frac{(1-|x|^{2})}{2}U(|x|e_{n}),
 \eeq where $$\Phi^{\ast}(x)=\Phi(x)-
\delta_{1}(|x|)K[\varphi_{1}](0)-\delta_{2}(|x|)H[\varphi_{2}](0).$$
If $x$ is not on the ray $[0,e_{n}]$, then we choose a unitary
transformation $\mathcal{O}$ such that $\mathcal{O}(e_{n})=x/|x|$.
By making use of the biharmonic mapping
$\mathcal{B}(y)=\Phi(\mathcal{O}(y))$, we get

\beqq\mathcal{B}(|x|e_{n})=\Phi(\mathcal{O}(|x|e_{n}))=\Phi(x),\eeqq
which, together with (\ref{eq-xj-8}), implies that, for any
$x\in\mathbb{B}^{n}$,

\beq\label{eq-xj-10z}
|\Phi^{\ast}(x)|&\leq&\frac{|n-4|}{4}\|\varphi_{1}\|_{\infty}(1-|x|^{2})U(|x|e_{n})
\\
\nonumber
&&+\frac{n}{4}\|\varphi_{1}\|_{\infty}U^{\ast}(|x|e_{n})+\|\varphi_{2}\|_{\infty}\frac{(1-|x|^{2})}{2}U(|x|e_{n}).
 \eeq

Next, we estimate $|G[g]|$.

\noindent $\mathbf{Case~1.}$ $n\neq2,4$.

For $x,y\in\mathbb{B}^{n}$, let
$z=\phi_{x}(y)\in\Aut(\mathbb{B}^{n})$. Then
\be\label{eq-j-1}|x-\phi_{x}(z)|
=\left|\frac{(z-x|z|^{2})(1-|x|^{2})}{[x,z]^{2}}\right|=\frac{|z|(1-|x|^{2})}{[x,z]},\ee
which, together with  (\ref{II}), gives

\beqq
G_{2,n}(x,y)&=&c_{n}\frac{\left(1-|\phi_{x}(y)|^{n-4}\right)}{|x-y|^{n-4}}-\frac{c_{n}(n-4)}{2}\frac{|\phi_{x}(y)|^{n-4}(1-|\phi_{x}(y)|^{2})}{|x-y|^{n-4}}\\
\nonumber&=&c_{n}\frac{\left(1-|z|^{n-4}\right)[x,z]^{n-4}}{|z|^{n-4}(1-|x|^{2})^{n-4}}-\frac{c_{n}(n-4)}{2}\frac{\left(1-|z|^{2}\right)[x,z]^{n-4}}{(1-|x|^{2})^{n-4}}.
\eeqq 
 By \cite[Lemma
2.1]{GS}, we know that $G_{2,n}(x,y)>0$.
 Then, by changing variables, we obtain

\be\label{eq-1ef}
\int_{\mathbb{B}^{n}}|G_{2,n}(x,y)|dV(y)=\int_{\mathbb{B}^{n}}G_{2,n}(x,y)dV(y)=\mathcal{G}_{1}-\mathcal{G}_{2},
\ee where
$$\mathcal{G}_{1}(x)=c_{n}(1-|x|^{2})^{4}\int_{\mathbb{B}^{n}}\frac{\left(1-|z|^{n-4}\right)}{|z|^{n-4}[x,z]^{n+4}}dV(z)$$
and
$$\mathcal{G}_{2}(x)=c_{n}\frac{(n-4)}{2}(1-|x|^{2})^{4}\int_{\mathbb{B}^{n}}\frac{\left(1-|z|^{2}\right)}{[x,z]^{n+4}}dV(z).$$

Now we estimate $\mathcal{G}_{1}$ and $\mathcal{G}_{2}$. Using the
spherical coordinates and Lemma \Ref{pro-1}, we obtain
\beq\label{eq-2ef}
\int_{\mathbb{S}^{n-1}}\frac{d\sigma(\zeta)}{|\rho
x-\zeta|^{4+n}}&=&\frac{\Gamma\big(\frac{n}{2}\big)}{\sqrt{\pi}\Gamma\big(\frac{n-1}{2}\big)}\int_{0}^{\pi}
\frac{\sin^{n-2}t}{\left(1+\rho^{2}|x|^{2}-2\rho|x|\cos
t\right)^{\frac{n+4}{2}}}dt\\ \nonumber
&=&\frac{\Gamma\big(\frac{n}{2}\big)}{\sqrt{\pi}\Gamma\big(\frac{n-1}{2}\big)}\cdot
\frac{\sqrt{\pi}\Gamma\big(\frac{n-1}{2}\big)}{\Gamma\big(\frac{n}{2}\big)}{_{2}}F_{1}\Big(\frac{n+4}{2},3;\frac{n}{2};\rho^{2}|x|^{2}\Big)\\
\nonumber &=& {_{2}}F_{1}\Big(\frac{n+4}{2},3;\frac{n}{2};\rho^{2}|x|^{2}\Big)\\
\nonumber
&=&\sum_{k=0}^{\infty}\frac{(k+1)(k+2)(n+2k)(n+2k+2)}{2n(n+2)}\rho^{2k}|x|^{2k},
\eeq which gives that

\beq\label{eq-3ef}
\mathcal{G}_{1}(x)&=&c_{n}A_{n-1}(1-|x|^{2})^{4}\int_{0}^{1}\rho^{3}(1-\rho^{n-4})
\left(\int_{\mathbb{S}^{n-1}}\frac{d\sigma(\zeta)}{|x\rho-\zeta|^{n+4}}\right)d\rho\\
\nonumber
&=&\frac{(n-4)c_{n}A_{n-1}(1-|x|^{2})^{4}}{4n(n+2)}\sum_{k=0}^{\infty}(k+1)(2k+2+n)|x|^{2k}\\
\nonumber
&=&\frac{(n-4)c_{n}A_{n-1}(1-|x|^{2})\big(n+2-(n-2)|x|^{2}\big)}{4n(n+2)}
\eeq and

\beq\label{eq-4ef}
\mathcal{G}_{2}(x)&=&\frac{(n-4)c_{n}A_{n-1}}{2}(1-|x|^{2})^{4}\int_{0}^{1}
\left(\int_{\mathbb{S}^{n-1}}\frac{\rho^{n-1}(1-\rho^{2})d\sigma(\zeta)}{|x\rho-\zeta|^{n+4}}\right)d\rho\\
\nonumber
&=&\frac{(n-4)c_{n}A_{n-1}(1-|x|^{2})^{4}}{2n(n+2)}\sum_{k=0}^{\infty}(k+1)(k+2)|x|^{2k}\\
\nonumber &=&\frac{(n-4)c_{n}A_{n-1}(1-|x|^{2})}{n(n+2)}. \eeq It
follows from (\ref{eq-1ef}), (\ref{eq-3ef}) and (\ref{eq-4ef}) that

\beq\label{Green-1}
|G[g](x)|&\leq&\|g\|_{\infty}\int_{\mathbb{B}^{n}}|G_{2,n}(x,y)|dV(y)=\|g\|_{\infty}\big(\mathcal{G}_{1}-\mathcal{G}_{2}\big)\\
\nonumber
&\leq&\|g\|_{\infty}\frac{(n-4)(n-2)c_{n}A_{n-1}}{4n(n+2)}(1-|x|^{2})^{2}\\
\nonumber &=&\frac{\|g\|_{\infty}}{8n(n+2)}(1-|x|^{2})^{2}.\eeq

\noindent $\mathbf{Case~2.}$ $n=2$.

For $x,y\in\mathbb{B}^{2}$, let
$z=\phi_{x}(y)\in\Aut(\mathbb{B}^{2})$. By (\ref{II}), (\ref{III})
and Lemma \Ref{pro-1}, we have

\beqq
\mathcal{G}_{3}(x)&=&\int_{\mathbb{B}^{2}}|x-y|^{2}\log\frac{|x-y|^{2}}{[x,y]^{2}}dV(y)=(1-|x|^{2})^{4}\int_{\mathbb{B}^{2}}\frac{|z|^{2}\log|z|^{2}}{[x,z]^{6}}dV(z)\\
&=&4\pi(1-|x|^{2})^{4}\int_{0}^{1}\rho^{3}\log\rho\left(\int_{\mathbb{S}^{1}}\frac{d\sigma(\zeta)}{|x\rho-\zeta|^{6}}\right)d\rho\\
&=&-4\pi(1-|x|^{2})^{4}\sum_{k=0}^{\infty}(k+1)^{2}(k+2)^{2}|x|^{2k}\int_{0}^{1}\rho^{3+2k}\log\rho
d\rho\\
&=&-\frac{A_{1}}{8}(1-|x|^{2})(1+|x|^{2}), \eeqq which, together
with $G_{2,2}(x,y)>0,$ yields that

\beq\label{Green-2}
|G[g](x)|&\leq&\|g\|_{\infty}\int_{\mathbb{B}^{2}}|G_{2,2}(x,y)|dV(y)=\|g\|_{\infty}\int_{\mathbb{B}^{2}}G_{2,2}(x,y)dV(y)\\
\nonumber
&=&\|g\|_{\infty}c_{2}\left(\mathcal{G}_{3}(x)+\int_{\mathbb{B}^{2}}(x\otimes
y)dV(y)\right)\\ \nonumber
&=&\|g\|_{\infty}c_{2}\left(\mathcal{G}_{3}(x)+\frac{A_{1}(1-|x|^{2})}{4}\right)=\|g\|_{\infty}\frac{(1-|x|^{2})^{2}}{64}.
\eeq

\noindent $\mathbf{Case~3.}$ $n=4$.

Since $G_{2,4}(x,y)>0~\mbox{for}~x,y\in\mathbb{B}^{4}$, by
(\ref{eq-2ef}), we see that

\beq\label{Green-3}
|G[g](x)|&\leq&\|g\|_{\infty}\int_{\mathbb{B}^{4}}|G_{2,4}(x,y)|dV(y)=\|g\|_{\infty}\int_{\mathbb{B}^{4}}G_{2,4}(x,y)dV(y)\\
\nonumber
&=&\|g\|_{\infty}c_{4}\int_{\mathbb{B}^{4}}\left(1-|z|^{2}+\log|z|^{2}\right)\frac{(1-|x|^{2})^{4}}{[x,z]^{8}}dV(z)\\
\nonumber
&=&A_{3}\|g\|_{\infty}c_{4}\int_{0}^{1}\rho^{3}\left(1-\rho^{2}+\log\rho^{2}\right)
\left(\int_{\mathbb{S}^{3}}\frac{(1-|x|^{2})^{4}d\sigma(\zeta)}{|x\rho-\zeta|^{8}}\right)d\rho\\
\nonumber
&=&-\frac{A_{3}\|g\|_{\infty}c_{4}(1-|x|^{2})^{4}}{24}\sum_{k=0}^{\infty}(k+1)|x|^{2k}\\
\nonumber &=&\frac{(1-|x|^{2})^{2}}{192}\|g\|_{\infty}, \eeq where
$z=\phi_{x}(y)$. Hence (\ref{eq-th-3}) follows from
(\ref{eq-xj-10z}), (\ref{Green-1}), (\ref{Green-2}) and
(\ref{Green-3}).

In particular, for $n=2,3,4$, we compute the values of
$U^{\ast}(re_{n})$ and $U(re_{n})$, repectively, where $r=|x|$.
 Let
$\zeta=(\zeta_{1},\ldots,\zeta_{n})\in\mathbb{S}^{n-1}$ such that
$\zeta_{n}=\cos\theta$, where $\theta$ is the angle between the
vector $x$ and $x_{n}$ axis. Let $m(r)=2r/(1+r^{2})$. Elementary
calculations lead to

\beqq\label{pp-1}
&&\frac{(1-r^{2})^{3}\sin^{n-2}\theta}{(1+r^{2}-2r\cos\theta)^{\frac{n+2}{2}}}-
\frac{(1-r^{2})^{3}\sin^{n-2}\theta}{(1+r^{2}+2r\cos\theta)^{\frac{n+2}{2}}}\\
\nonumber
&=&\frac{(1-r^{2})^{3}}{(1+r^{2})^{\frac{n+2}{2}}}\sum_{k=0}^{\infty}\left(-\frac{(2+n)}{2}\right)_{k}\big((-1)^{k}-1\big)m^{k}(r)\cos^{k}\theta\sin^{n-2}\theta.
\eeqq and $$\int_{0}^{\frac{\pi}{2}}\cos^{k}\theta\sin^{n-2}\theta
d\theta=\frac{\Gamma\left(\frac{1+k}{2}\right)\Gamma\left(\frac{n-1}{2}\right)}{2\Gamma\left(\frac{n+k}{2}\right)},$$
which imply that

\beqq
U^{\ast}(re_{n})&=&\frac{\Gamma\left(\frac{n}{2}\right)}{\sqrt{\pi}\Gamma\left(\frac{n-1}{2}\right)}
\int_{0}^{\pi}\frac{(1-r^{2})^{3}\sin^{n-2}\theta}{(1+r^{2}-2r\cos\theta)^{\frac{n+2}{2}}}(\chi_{\mathbb{S}^{n-1}_{+}}-\chi_{\mathbb{S}^{n-1}_{-}})d\theta\\
&=&\frac{\Gamma\left(\frac{n}{2}\right)}{\sqrt{\pi}\Gamma\left(\frac{n-1}{2}\right)}
\int_{0}^{\frac{\pi}{2}}\left(\frac{(1-r^{2})^{3}\sin^{n-2}\theta}{(1+r^{2}-2r\cos\theta)^{\frac{n+2}{2}}}-
\frac{(1-r^{2})^{3}\cos^{n-2}\theta}{(1+r^{2}+2r\sin\theta)^{\frac{n+2}{2}}}\right)d\theta\\
&=&\frac{\Gamma\left(\frac{n}{2}\right)}{\sqrt{\pi}\Gamma\left(\frac{n-1}{2}\right)}
\int_{0}^{\frac{\pi}{2}}\left(\frac{(1-r^{2})^{3}\sin^{n-2}\theta}{(1+r^{2}-2r\cos\theta)^{\frac{n+2}{2}}}-
\frac{(1-r^{2})^{3}\sin^{n-2}\theta}{(1+r^{2}+2r\cos\theta)^{\frac{n+2}{2}}}\right)d\theta\\
&=&Q(r), \eeqq where
$$Q(r)=\frac{\Gamma\left(\frac{n}{2}\right)}{\sqrt{\pi}\Gamma\left(\frac{n-1}{2}\right)}\frac{(1-r^{2})^{3}}{(1+r^{2})^{\frac{n+2}{2}}}
\sum_{k=0}^{\infty}\frac{\Gamma\left(\frac{1+k}{2}\right)\Gamma\left(\frac{n-1}{2}\right)}{2\Gamma\left(\frac{n+k}{2}\right)}
\left(-\frac{(2+n)}{2}\right)_{k}\big((-1)^{k}-1\big)m^{k}(r).$$
Since $(-1)^{k}-1=0$ for $k=2j$, where $j\in\{0,1,\ldots,\}$, we see
that $Q(r)$ can be rewritten as

\beq\label{vop-1}\nonumber
Q(r)&=&-\frac{\Gamma\left(\frac{n}{2}\right)}{\sqrt{\pi}}(1-r^{2})^{3}\sum_{j=0}^{\infty}\frac{(2r)^{2j+1}\Gamma(1+j)}{\Gamma\left(\frac{n+2j+1}{2}\right)(1+r^{2})^{\frac{n}{2}+2j+2}}
\left(-\frac{(2+n)}{2}\right)_{2j+1}\\
&=&\frac{\Gamma\left(\frac{n}{2}\right)2(2+n)r(1-r^{2})^{3}{_{4}}F_{3}\left(\big\{1,1+\frac{n}{4},\frac{3}{2}+\frac{n}{4}\big\},
\big\{\frac{3}{2},\frac{1}{2}+\frac{n}{2}\big\},
\frac{4r^{2}}{(1+r^{2})^{2}}\right)}{\sqrt{\pi}\Gamma\left(\frac{n-1}{2}\right)(n-1)(1+r^{2})^{2+\frac{n}{2}}},
\eeq where ${_{4}}F_{3}$ is defined in \cite{GR}.
 By applying \cite[Eq. 3.1.8]{GR} to (\ref{vop-1}), we obtain  the values of $U^{\ast}(re_{n})$ (see the Table 1). The values of $U(re_{n})$ follows from \cite[Remark 2.7]{K5}
 (see also the Table 1).
The proof of this theorem is complete. \qed







\section{ modulus of continuity  of solutions to the inhomogeneous
 biharmonic Dirichlet problems}\label{csw-sec4}

We begin this part with the following two Lemmas which will be used
in the proof of Theorem \ref{thm-2}.

\begin{lem}\label{c1-1}
For $x\in\mathbb{B}^{n}$,
$$\int_{\mathbb{S}^{n-1}}K_{n}(x,\zeta)d\sigma(\zeta)=1.$$
\end{lem}
\bpf By the spherical coordinate transformation (see section
\ref{sbcsw-sec2.4}) and Lemma \Ref{pro-1}, we have

\beq\label{c-20v}
\int_{\mathbb{S}^{n-1}}\frac{d\sigma(\zeta)}{|x-\zeta|^{2+n}}&=&\frac{1}{\int_{0}^{\pi}\sin^{n-2}
t}\int_{0}^{\pi}\frac{\sin^{n-2} t\, dt}{\big(1+|x|^{2}-2|x|\cos
t\big)^{\frac{n+2}{2}}}\\ \nonumber &=&
\sum_{j=0}^{\infty}\frac{(n+2k)(k+1)}{n}|x|^{2k}\\
\nonumber
&=&\frac{1}{n}\left(\frac{n}{(1-|x|^{2})^{2}}+\frac{4|x|^{2}}{(1-|x|^{2})^{3}}\right).\eeq
Elementary computations show that
$$\int_{\mathbb{S}^{n-1}}\frac{(1-|x|^{2})^{2}}{|x-\zeta|^{n}}d\sigma(\zeta)=1-|x|^{2},$$
which, together with (\ref{c-20v}), implies that

\beqq\int_{\mathbb{S}^{n-1}}K_{n}(x,\zeta)d\sigma(\zeta)&=&\frac{n}{4}\int_{\mathbb{S}^{n-1}}\frac{(1-|x|^{2})^{3}}{|x-\zeta|^{n+2}}d\sigma(\zeta)-
\frac{(n-4)}{4}\int_{\mathbb{S}^{n-1}}\frac{(1-|x|^{2})^{2}}{|x-\zeta|^{n}}d\sigma(\zeta)\\
&=&\frac{(1-|x|^{2})^{3}}{4}\left(\frac{n}{(1-|x|^{2})^{2}}+\frac{4|x|^{2}}{(1-|x|^{2})^{3}}\right)-\frac{(n-4)}{4}(1-|x|^{2})\\
&=&1. \eeqq The proof of this lemma is finished. \epf

\begin{Lem}\label{K-3}{\rm (\cite[Lemma 2.5]{K-2011})}
Let $\varrho$ be a bounded (absolutely) integrable function defined
on a bounded domain $\Omega\subset\mathbb{R}^{n}$. Then the
potential type integral

$$\tau(x)=\int_{\Omega}\frac{\varrho(y)dV(y)}{|x-y|^{\alpha}}$$ belongs to
the space $\mathcal{C}^{k}(\mathbb{R}^{n},\mathbb{R})$, where
$k+\alpha<n.$ Moreover,

$$\nabla\tau(x)=\int_{\Omega}\nabla\left(\frac{1}{|x-y|^{\alpha}}\varrho(y)\right)dV(y).$$
\end{Lem}

\subsection*{The proof of Theorem \ref{thm-2}} We divide the proof of this theorem into four
steps. 

\bst\label{bst-1.0} The estimate of $|D_{K[\varphi_{1}]}|$. \est


For $x=(x_{1},\ldots,x_{n})\in\mathbb{B}^{n}$ and
$\zeta=(\zeta_{1},\ldots,\zeta_{n})\in\mathbb{S}^{n-1}$, we obtain

\beqq \frac{\partial}{\partial
x_{k}}K_{n}(x,\zeta)&=&\frac{1}{4}\bigg(-\frac{6n(1-|x|^{2})^{2}x_{k}}{|x-\zeta|^{n+2}}-\frac{n(n+2)(1-|x|^{2})^{3}(x_{k}-\zeta_{k})}{|x-\zeta|^{n+4}}\\
\nonumber
&&+\frac{4(n-4)(1-|x|^{2})x_{k}}{|x-\zeta|^{n}}+\frac{n(n-4)(1-|x|^{2})^{2}(x_{k}-\zeta_{k})}{|x-\zeta|^{n+2}}\bigg),
\eeqq where $k\in\{1,\ldots,n\}$.

Then, for any
$\xi=(\xi_{1},\ldots,\xi_{n})\in\mathbb{R}^{n}\setminus\{0\}$, we
obtain

\beqq\langle\nabla K_{n}(x,\zeta),\xi\rangle&=&\frac{1}{4}\sum_{k=1}^{n}\bigg(-\frac{6n(1-|x|^{2})^{2}x_{k}\xi_{k}}{|x-\zeta|^{n+2}}-
\frac{n(n+2)(1-|x|^{2})^{3}(x_{k}-\zeta_{k})\xi_{k}}{|x-\zeta|^{n+4}}\\
\nonumber
&&+\frac{4(n-4)(1-|x|^{2})x_{k}\xi_{k}}{|x-\zeta|^{n}}+\frac{n(n-4)(1-|x|^{2})^{2}(x_{k}-\zeta_{k})\xi_{k}}{|x-\zeta|^{n+2}}\bigg),
 \eeqq
which, together with Cauchy-Schwarz's inequality, implies that

\beq\label{vp-1} |\langle\nabla
K_{n}(x,\zeta),\xi\rangle|&\leq&\frac{3n}{2}\frac{(1-|x|^{2})^{2}|x||\xi|}{|x-\zeta|^{n+2}}+
\frac{n(n+2)}{4}\frac{(1-|x|^{2})^{3}|\xi|}{|x-\zeta|^{n+3}}\\
\nonumber
&&+\frac{|n-4|(1-|x|^{2})|x||\xi|}{|x-\zeta|^{n}}+\frac{n|n-4|(1-|x|^{2})^{2}|\xi|}{4|x-\zeta|^{n+1}},
 \eeq
where $\langle\cdot,\cdot\rangle$ is the Euclidean inner product.

By (\ref{vp-1}) and Lemma \ref{c1-1}, for
$x=|x|x^{\ast}\in\mathbb{B}^{n}\setminus\{0\}$, we get

\beq\label{c-21cg}|D_{K[\varphi_{1}]}(x)\xi|&=&\left|\int_{\mathbb{S}^{n-1}}\langle\nabla
K_{n}(x,\zeta),\xi\rangle\varphi_{1}(\zeta)d\sigma(\zeta)\right|\\
\nonumber &=&\left|\int_{\mathbb{S}^{n-1}}\langle\nabla
K_{n}(x,\zeta),\xi\rangle(\varphi_{1}(\zeta)-\varphi_{1}(x^{\ast}))d\sigma(\zeta)\right|\\
\nonumber &\leq&|\xi|\sum_{j=1}^{4}I_{j}(x), \eeq where
$$I_{1}(x)=\frac{3n}{2}\int_{\mathbb{S}^{n-1}}\frac{(1-|x|^{2})^{2}|x|}{|x-\zeta|^{n+2}}|\varphi_{1}(\zeta)-\varphi_{1}(x^{\ast})|d\sigma(\zeta),$$
$$I_{2}(x)=\frac{n(n+2)}{4}\int_{\mathbb{S}^{n-1}}\frac{(1-|x|^{2})^{3}}{|x-\zeta|^{n+3}}|\varphi_{1}(\zeta)-\varphi_{1}(x^{\ast})|d\sigma(\zeta),$$
$$I_{3}(x)=|n-4|\int_{\mathbb{S}^{n-1}}\frac{(1-|x|^{2})|x|}{|x-\zeta|^{n}}|\varphi_{1}(\zeta)-\varphi_{1}(x^{\ast})|d\sigma(\zeta)$$
and
$$I_{4}(x)=\frac{n|n-4|}{4}\int_{\mathbb{S}^{n-1}}\frac{(1-|x|^{2})^{2}}{|x-\zeta|^{n+1}}|\varphi_{1}(\zeta)-\varphi_{1}(x^{\ast})|d\sigma(\zeta).$$

Now we first estimate $I_{1}$ and $I_{2}$. By (\ref{c-20v}) and
Cauchy-Schwarz's inequality, we have

\beq\label{eq-ca-k0} \nonumber
\int_{\mathbb{S}^{n-1}}\frac{(1-|x|^{2})^{2}}{|x-\zeta|^{n+1}}d\sigma(\zeta)&\leq&\left(\int_{\mathbb{S}^{n-1}}\frac{1-|x|^{2}}{|x-\zeta|^{n}}d\sigma(\zeta)\right)^{\frac{1}{2}}
\left(\int_{\mathbb{S}^{n-1}}\frac{(1-|x|^{2})^{3}}{|x-\zeta|^{n+2}}d\sigma(\zeta)\right)^{\frac{1}{2}}\\
 &=&\left(1-|x|^{2}+\frac{4}{n}|x|^{2}\right)^{\frac{1}{2}}. \eeq

 Since
$\varphi_{1}\in
\mathcal{L}_{\omega}(\mathbb{S}^{n-1},\mathbb{R}^{m})$, we see that
there is positive constant $L$ such that
\be\label{eq-ca-k1}|\varphi_{1}(\zeta)-\varphi_{1}(x^{\ast})|\leq
L\omega(|\zeta-x^{\ast}|),~\zeta\in\mathbb{S}^{n-1}.\ee Simple
calculations show that, for $\zeta\in\mathbb{S}^{n-1}$,
\be\label{eq-ca-k2}|\zeta-x^{\ast}|\leq|\zeta-x|+|x-x^{\ast}|=|\zeta-x|+(1-|x|)\leq2|\zeta-x|.\ee
By (\ref{eq-ca-k0}), (\ref{eq-ca-k1}) and (\ref{eq-ca-k2}), we have

\beq\label{eq-ca-k4}
I_{1}(x)&=&\frac{3n|x|}{2}\int_{\mathbb{S}^{n-1}}\frac{(1-|x|^{2})^{2}}{|x-\zeta|^{n+2}}|\varphi_{1}(\zeta)-\varphi_{1}(x^{\ast})|d\sigma(\zeta)\\
\nonumber
&\leq&\frac{3n|x|L}{2}\int_{\mathbb{S}^{n-1}}\frac{(1-|x|^{2})^{2}}{|x-\zeta|^{n+2}}\frac{\omega(|\zeta-x^{\ast}|)}{|\zeta-x^{\ast}|}|\zeta-x^{\ast}|d\sigma(\zeta)
\\
\nonumber &\leq&
3n|x|Lc\int_{\mathbb{S}^{n-1}}\frac{(1-|x|^{2})^{2}}{|x-\zeta|^{n+1}}d\sigma(\zeta)\\
\nonumber
&\leq&3nLc|x|\left(1-|x|^{2}+\frac{4}{n}|x|^{2}\right)^{\frac{1}{2}}=M_{1}(n),
 \eeq
where
$$M_{1}(n)=
\begin{cases}
\displaystyle  6\sqrt{n}Lc, & 2\leq n\leq 4,\\
\displaystyle \\ \displaystyle nLc\sqrt{\frac{n}{(n-4)}},& n>4.
\end{cases}$$

Applying (\ref{c-20v}), (\ref{eq-ca-k1}) and (\ref{eq-ca-k2}), we
obtain \beq\label{eq-ca-k5}
I_{2}(x)&\leq&\frac{n(n+2)L}{4}\int_{\mathbb{S}^{n-1}}\frac{(1-|x|^{2})^{3}}{|x-\zeta|^{n+3}}\frac{\omega(|\zeta-x^{\ast}|)}{|\zeta-x^{\ast}|}|\zeta-x^{\ast}|d\sigma(\zeta)\\
\nonumber&\leq&\frac{n(n+2)Lc}{2}\int_{\mathbb{S}^{n-1}}\frac{(1-|x|^{2})^{3}}{|x-\zeta|^{n+2}}d\sigma(\zeta)\\
\nonumber&=&\frac{n(n+2)Lc}{2}\left(1-|x|^{2}+\frac{4}{n}|x|^{2}\right)=M_{2}(n),
 \eeq
where
$$M_{2}(n)=
\begin{cases}
\displaystyle 2(n+2)Lc, & 2\leq n\leq 4,\\
\displaystyle \\ \displaystyle \frac{n(n+2)Lc}{2},& n>4.
\end{cases}$$

Next, we estimate $I_{3}$ and $I_{4}$. By (\ref{eq-ca-k1}) and
(\ref{eq-ca-k2}), we get

\be\label{eq-ca-k6}I_{3}(x)\leq|n-4|L|x|\int_{\mathbb{S}^{n-1}}\frac{(1-|x|^{2})}{|x-\zeta|^{n}}\omega(|\zeta-x^{\ast}|)d\sigma(\zeta)\leq|n-4|L\omega(2)\ee
and
\beq\label{eq-ca-k7}I_{4}(x)&\leq&\frac{n|n-4|L}{4}\int_{\mathbb{S}^{n-1}}\frac{(1-|x|^{2})^{2}}{|x-\zeta|^{n+1}}\frac{\omega(|\zeta-x^{\ast}|)}{|\zeta-x^{\ast}|}|\zeta-x^{\ast}|d\sigma(\zeta)\\
\nonumber&\leq&\frac{n|n-4|Lc}{4}\int_{\mathbb{S}^{n-1}}\frac{(1-|x|^{2})^{2}}{|x-\zeta|^{n+1}}|\zeta-x^{\ast}|d\sigma(\zeta)
\\
\nonumber&\leq&\frac{n|n-4|Lc}{2}\int_{\mathbb{S}^{n-1}}\frac{(1-|x|^{2})^{2}}{|x-\zeta|^{n}}d\sigma(\zeta)
\\
\nonumber&\leq&\frac{n|n-4|Lc}{2}.
 \eeq
It follows from (\ref{c-21cg}),  (\ref{eq-ca-k4}), (\ref{eq-ca-k5}),
(\ref{eq-ca-k6}) and (\ref{eq-ca-k7}) that

\be\label{eq-ca-k8}|D_{K[\varphi_{1}]}(x)\xi|\leq
\left(M_{1}(n)+M_{2}(n)+|n-4|L\omega(2)+\frac{n|n-4|Lc}{2}\right)|\xi|.\ee

\bst\label{bst-2.0} The estimate  of $|D_{H[\varphi_{2}]}|$. \est
For $x=(x_{1},\ldots,x_{n})\in\mathbb{B}^{n}$ and
$\zeta=(\zeta_{1},\ldots,\zeta_{n})\in\mathbb{S}^{n-1}$, we obtain
\beqq \frac{\partial}{\partial
x_{k}}H_{n}(x,\zeta)=-\frac{1}{2}\left(\frac{4(1-|x|^{2})x_{k}}{|x-\zeta|^{n}}+\frac{n(1-|x|^{2})^{2}(x_{k}-\zeta_{k})}{|x-\zeta|^{n+2}}\right),
\eeqq where $k\in\{1,\ldots,n\}$. Then, for any
$\xi=(\xi_{1},\ldots,\xi_{n})\in\mathbb{R}^{n}\setminus\{0\}$, we
have

\beqq \langle\nabla
H_{n}(x,\zeta),\xi\rangle=-\frac{1}{2}\sum_{k=1}^{n}\left(\frac{4(1-|x|^{2})x_{k}\xi_{k}}{|x-\zeta|^{n}}+\frac{n(1-|x|^{2})^{2}(x_{k}-\zeta_{k})\xi_{k}}{|x-\zeta|^{n+2}}\right),
\eeqq which, together with Cauchy-Schwarz's inequality, gives that

\beq\label{eq-ca-k10} |\langle\nabla
H_{n}(x,\zeta),\xi\rangle|\leq\left(\frac{2(1-|x|^{2})|x|}{|x-\zeta|^{n}}+\frac{n}{2}\frac{(1-|x|^{2})^{2}}{|x-\zeta|^{n+1}}\right)|\xi|.
\eeq Applying (\ref{eq-ca-k0}) and (\ref{eq-ca-k10}), we see that

\beq\label{eq-c-h10}|D_{H[\varphi_{2}]}(x)\xi|&=&\left|\int_{\mathbb{S}^{n-1}}\langle\nabla
H_{n}(x,\zeta),\xi\rangle\varphi_{2}(\zeta)d\sigma(\zeta)\right|\\
\nonumber &\leq&\left(2|x|+
\frac{n}{2}\int_{\mathbb{S}^{n-1}}\frac{(1-|x|^{2})^{2}}{|x-\zeta|^{n+1}}d\sigma(\zeta)\right)|\xi|\|\varphi_{2}\|_{\infty}\\
\nonumber&\leq& \left(2+
\frac{n}{2}\left(1-|x|^{2}+\frac{4}{n}|x|^{2}\right)^{\frac{1}{2}}\right)|\xi|\|\varphi_{2}\|_{\infty}\\
\nonumber &=&M_{3}(n)|\xi|\|\varphi_{2}\|_{\infty},\eeq where
$$M_{3}(n)=
\begin{cases}
\displaystyle 2+\sqrt{n}, & 2\leq n\leq 4,\\
\displaystyle \\ \displaystyle \frac{4+n}{2},& n>4.
\end{cases}$$

\bst\label{bst-3.0} The estimate  of $|D_{G[g]}|$. \est

\noindent $\mathbf{Case~1.}$ $n\neq2,4$.

For $x=(x_{1},\ldots,x_{n}),y=(y_{1},\ldots,y_{n})\in\mathbb{B}^{n}$
with $x\neq y$,  we have

\beqq \frac{\partial}{\partial
x_{k}}G_{2,n}(x,y)&=&c_{n}(4-n)\bigg(|x-y|^{2-n}(x_{k}-y_{k})-[x,y]^{2-n}\left(x_{k}|y|^{2}-y_{k}\right)\\
&&-(1-|y|^{2})[x,y]^{2-n}x_{k}+\frac{(2-n)}{2}[x,y]^{-n}(x\otimes
y)\left(x_{k}|y|^{2}-y_{k}\right)\bigg). \eeqq Then, for any
$\xi=(\xi_{1},\ldots,\xi_{n})\in\mathbb{R}^{n}\setminus\{0\}$, we
get

\beq\label{eq-c-h11} |\langle\nabla
G_{2,n}(x,y),\xi\rangle|&\leq&|c_{n}||4-n|\bigg(|x-y|^{3-n}+(1-|y|^{2})[x,y]^{2-n}|x|\\
\nonumber &&+[x,y]^{3-n}|y|+\frac{(n-2)}{2}[x,y]^{1-n}(x\otimes
y)|y|\bigg)|\xi|, \eeq which, together with Lemma \Ref{K-3}, yields
that

\beq\label{eq-ca-h13}|D_{G[g]}(x)\xi|&\leq&\|g\|_{\infty}\int_{\mathbb{B}^{n}}\left|\langle\nabla
G_{2,n}(x,\zeta),\xi\rangle\right| dV(y)\\ \nonumber
&\leq&\|g\|_{\infty}|c_{n}||4-n||\xi|\sum_{j=5}^{8}I_{j}(x),
 \eeq
where
$$I_{5}(x)=\int_{\mathbb{B}^{n}}\frac{dV(y)}{|x-y|^{n-3}},~I_{6}(x)=\int_{\mathbb{B}^{n}}\frac{dV(y)}{[x,y]^{n-3}},~I_{7}(x)=|x|\int_{\mathbb{B}^{n}}\frac{(1-|y|^{2})}{[x,y]^{n-2}}dV(y),$$
and $$I_{8}(x)=\frac{(n-2)}{2}\int_{\mathbb{B}^{n}}\frac{x\otimes
y}{[x,y]^{n-1}}dV(y).$$

Now we estimate $I_{5}$. Let $z=\phi_{x}(y)\in\Aut(\mathbb{B}^{n})$.
Then, by (\ref{III}) and (\ref{eq-j-1}), we have  

\beq\label{eq-ca-h16}
I_{5}(x)&=&\int_{\mathbb{B}^{n}}\frac{|x-y|}{|x-y|^{n-2}}dV(y)\leq2\int_{\mathbb{B}^{n}}\frac{dV(y)}{|x-y|^{n-2}}\\
\nonumber
&=&2\int_{\mathbb{B}^{n}}\frac{1}{|x-\phi_{x}(z)|^{n-2}}\frac{(1-|x|^{2})^{n}}{[x,z]^{2n}}dV(z)\\
\nonumber
&=&2(1-|x|^{2})^{2}\int_{\mathbb{B}^{n}}\frac{dV(z)}{|z|^{n-2}[x,z]^{n+2}}.
\eeq It follows from (\ref{c-20v}) that

\beqq\label{eq-ca-h17}
\int_{\mathbb{B}^{n}}\frac{dV(z)}{|z|^{n-2}[x,z]^{n+2}}&=&A_{n-1}\int_{0}^{1}\left(\rho\int_{\mathbb{S}^{n-1}}\frac{d\sigma(\zeta)}{|x\rho-\zeta|^{n+2}}\right)d\rho\\
\nonumber
&=&A_{n-1}\int_{0}^{1}\left(\frac{\rho}{(1-|x|^{2}\rho^{2})^{2}}+\frac{4}{n}\frac{|x|^{2}\rho^{3}}{(1-|x|^{2}\rho^{2})^{3}}\right)d\rho\\
\nonumber
&=&\frac{A_{n-1}}{2}\frac{1}{(1-|x|^{2})}+\frac{A_{n-1}}{n}\frac{|x|^{2}}{(1-|x|^{2})^{2}},
\eeqq which, together with (\ref{eq-ca-h16}), implies that
\be\label{eq-ca-h18} I_{5}(x)\leq
\frac{A_{n-1}\left(n-(n-2)|x|^{2}\right)}{n}\leq A_{n-1}. \ee

Next, we estimate $I_{6}$. By (\ref{II}), we have

$$
|z|=\big|\phi_{x}(\phi_{x}(z))\big|=\frac{|x-\phi_{x}(z)|}{[x,\phi_{x}(z)]},
$$
which gives that

\be\label{eq-ca-h19}\frac{1}{[x,\phi_{x}(z)]}=\frac{|z|}{|x-\phi_{x}(z)|}.
\ee Applying (\ref{eq-j-1}) and (\ref{eq-ca-h19}), we see that

\beq\label{eq-nw-1}
\int_{\mathbb{B}^{n}}\frac{dV(y)}{[x,y]^{n-2}}&=&\int_{\mathbb{B}^{n}}\frac{1}{[x,\phi_{x}(z)]^{n-2}}\frac{(1-|x|^{2})^{n}}{[x,z]^{2n}}dV(z)\\
\nonumber
&=&(1-|x|^{2})^{2}\int_{\mathbb{B}^{n}}\frac{dV(z)}{[x,z]^{n+2}}\\
\nonumber
&\leq&(1-|x|^{2})^{2}\int_{\mathbb{B}^{n}}\frac{dV(z)}{|z|^{n-2}[x,z]^{n+2}}\leq
\frac{A_{n-1}}{2},\eeq which, together with (\ref{eq-j-1}) and
(\ref{eq-ca-h17}), implies that

\be\label{eq-ca-h21}
I_{6}(x)\leq2\int_{\mathbb{B}^{n}}\frac{dV(y)}{[x,y]^{n-2}}
 \leq A_{n-1}. \ee

At last, we estimate $I_{7}$ and $I_{8}$. It follows from
 (\ref{eq-nw-1})
that

\be\label{eq-ca-h22}
I_{7}(x)\leq\int_{\mathbb{B}^{n}}\frac{dV(y)}{[x,y]^{n-2}}\leq\frac{A_{n-1}}{2}.
 \ee

Since
$$\frac{x\otimes
y}{(1-|x||y|)}\leq\frac{x\otimes y}{(1-|x||y|)^{2}}\leq1,$$ by
(\ref{eq-nw-1}), we see that

\beq\label{eq-ca-h23}
I_{8}(x)&\leq&\frac{(n-2)}{2}\int_{\mathbb{B}^{n}}\frac{x\otimes
y}{[x,y]^{n-2}(1-|x||y|)}dV(y)\leq\frac{(n-2)}{2}\int_{\mathbb{B}^{n}}\frac{dV(y)}{[x,y]^{n-2}}\\
\nonumber &\leq&\frac{(n-2)A_{n-1}}{4}.
 \eeq

By (\ref{eq-ca-h18}),  (\ref{eq-ca-h21}), (\ref{eq-ca-h22}) and
(\ref{eq-ca-h23}), we conclude that

\be\label{eq-ca-g1}|D_{G[g]}(x)\xi|
\leq\|g\|_{\infty}|c_{n}||4-n|A_{n-1}\left(\frac{5}{2}+\frac{n-2}{4}\right)|\xi|.
 \ee

\noindent $\mathbf{Case~2.}$ $n=2$.

 For
$x=(x_{1},x_{2}),y=(y_{1},y_{2})\in\mathbb{B}^{n}$ with $x\neq y$,
 we have

\beqq \frac{\partial}{\partial
x_{k}}G_{2,2}(x,y)&=&2c_{2}\bigg((x_{k}-y_{k})\log\frac{|x-y|^{2}}{[x,y]^{2}}+(x_{k}-y_{k})\\
&&-\frac{|x-y|^{2}}{[x,y]^{2}}\big(x_{k}|y|^{2}-y_{k}\big)-x_{k}(1-|y|^{2})\bigg),
\eeqq which gives that,  for any
$\xi=(\xi_{1},\xi_{2})\in\mathbb{R}^{2}\setminus\{0\}$,

\beq\label{eq-hj-1} |\langle\nabla
G_{2,2}(x,y),\xi\rangle|&\leq&2c_{2}\bigg(|x-y|\bigg|\log\frac{|x-y|^{2}}{[x,y]^{2}}\bigg|+|x-y|\\
\nonumber &&+\frac{|x-y|^{2}|y|}{[x,y]}+|x|(1-|y|^{2})\bigg)|\xi|,
\eeq where $k\in\{1,2\}$.

Then, by (\ref{eq-hj-1}) and  Lemma \Ref{K-3}, we have

\beq\label{eq-hj-2}|D_{G[g]}(x)\xi|&\leq&\|g\|_{\infty}\int_{\mathbb{B}^{2}}\left|\langle\nabla
G_{2,2}(x,\zeta),\xi\rangle\right| dV(y)\\ \nonumber
&\leq&2c_{2}\|g\|_{\infty}\left(I_{9}(x)+2\pi+I_{10}(x)+\frac{\pi}{4}|x|\right)|\xi|,
 \eeq where $$I_{9}(x)=\int_{\mathbb{B}^{2}}|x-y|\bigg|\log\frac{|x-y|^{2}}{[x,y]^{2}}\bigg|dV(y)~\mbox{and}~I_{10}(x)=\int_{\mathbb{B}^{2}}\frac{|x-y|^{2}}{[x,y]}dV(y).$$

Now we estimate $I_{9}$. By Cauchy-Schwarz's inequality and Lemma
\Ref{pro-1}, we see that

\beq\label{eq-hj-4}\int_{\mathbb{S}^{1}}\frac{d\sigma(\zeta)}{|x\rho-\zeta|^{5}}&\leq&\left(\int_{\mathbb{S}^{1}}\frac{d\sigma(\zeta)}{|x\rho-\zeta|^{4}}\right)^{\frac{1}{2}}
\left(\int_{\mathbb{S}^{1}}\frac{d\sigma(\zeta)}{|x\rho-\zeta|^{6}}\right)^{\frac{1}{2}}\\
\nonumber
&\leq&\sum_{j=0}^{\infty}\frac{(n+1)^{2}(n+2)^{2}}{4}|x|^{2j}\rho^{2j}.
\eeq

Let $z=\phi_{x}(y)\in\Aut(\mathbb{B}^{n})$. Then, by
 (\ref{eq-j-1}) and (\ref{eq-hj-4}), we get

\beq\label{eq-hj-5}
I_{9}(x)&=&\int_{\mathbb{B}^{2}}|x-y|\log\frac{[x,y]^{2}}{|x-y|^{2}}dV(y)=\int_{\mathbb{B}^{2}}\frac{(1-|x|^{2})^{3}|z|}{[x,z]^{5}}\log\frac{1}{|z|^{2}}dV(z)\\
\nonumber
&=&2\pi(1-|x|^{2})^{3}\int_{0}^{1}\rho^{2}\log\frac{1}{\rho^{2}}\left(\int_{\mathbb{S}^{1}}\frac{d\sigma(\zeta)}{|x\rho-\zeta|^{5}}\right)d\rho\\
\nonumber
&\leq&2\pi(1-|x|^{2})^{3}\sum_{j=0}^{\infty}\frac{(n+1)^{2}(n+2)^{2}}{4}|x|^{2j}\int_{0}^{1}\rho^{2+2j}\log\frac{1}{\rho^{2}}d\rho\\
\nonumber
&\leq&2\pi(1-|x|^{2})^{3}\sum_{j=0}^{\infty}\frac{(n+1)(n+2)}{8}|x|^{2j}=\frac{\pi}{2}.
\eeq

At last, we estimate $I_{10}$. Since
$\phi_{x}(y)\in\Aut(\mathbb{B}^{n})$, we see that

\be\label{eq-hj-7}
I_{10}(x)=\int_{\mathbb{B}^{2}}|x-y||\phi_{x}(y)|dV(y)\leq2\int_{\mathbb{B}^{2}}dV(y)=2\pi.
\ee

Hence, in this case, it follows from (\ref{eq-hj-2}),
(\ref{eq-hj-5}) and (\ref{eq-hj-7}) that there is a positive
constant $M_{4}$ such that

\be\label{eq-hj-87}|D_{G[g]}(x)\xi|\leq M_{4}\|g\|_{\infty}|\xi|.\ee

\noindent $\mathbf{Case~3.}$ $n=4$.

For $x=(x_{1},x_{2}, x_{3}, x_{4} ),y=(y_{1},y_{2},y_{3},
y_{4})\in\mathbb{B}^{n}$ with $x\neq y$,  we have

\beqq \frac{\partial}{\partial
x_{k}}G_{2,4}(x,y)&=&2c_{4}\bigg(\frac{x_{k}-y_{k}}{|x-y|^{2}}-\frac{x_{k}|y|^{2}-y_{k}}{[x,y]^{2}}-\frac{(1-|y|^{2})x_{k}}{[x,y]^{2}}\\
&&-\frac{(x\otimes
y)\big(x_{k}|y|^{2}-y_{k}\big)}{[x,y]^{4}}\bigg)\\ &=&
2c_{4}\bigg(\frac{x_{k}(1-|y|^{2})\big(1+|y|^{2}-2\langle
x,y\rangle\big)-y_{k}(x\otimes
y)}{|x-y|^{2}[x,y]^{2}}\\
&&-\frac{(1-|y|^{2})x_{k}}{[x,y]^{2}}-\frac{(x\otimes
y)\big(x_{k}|y|^{2}-y_{k}\big)}{[x,y]^{4}}\bigg), \eeqq which yields
that, for any
$\xi=(\xi_{1},\xi_{2},\xi_{3},\xi_{4})\in\mathbb{R}^{n}\setminus\{0\}$,

\beq\label{eq-hj-s1} |\langle\nabla
G_{2,4}(x,y),\xi\rangle|&\leq&2|c_{4}|\bigg(\frac{(1-|y|^{2})\big(1+|y|^{2}-2\langle
x,y\rangle\big)|x|+(x\otimes y)|y|}{|x-y|^{2}[x,y]^{2}}\\ \nonumber
&&+\frac{(1-|y|^{2})|x|}{[x,y]^{2}}+\frac{(x\otimes
y)|y|}{[x,y]^{3}}\bigg)|\xi|\\ \nonumber
 &\leq&2|c_{4}|\bigg(\frac{(1-|y|^{2})\big(1+|y|^{2}-2\langle
x,y\rangle\big)+(x\otimes y)}{|x-y|^{2}[x,y]^{2}}\\ \nonumber
&&+\frac{(1-|y|^{2})}{[x,y]^{2}}+\frac{x\otimes
y}{[x,y]^{3}}\bigg)|\xi|\\ \nonumber
&=&2|c_{4}|\bigg(\frac{(1-|y|^{2})}{[x,y]^{2}}+\frac{2(x\otimes
y)}{|x-y|^{2}[x,y]^{2}}\\
\nonumber &&+\frac{(1-|y|^{2})}{[x,y]^{2}}+\frac{x\otimes
y}{[x,y]^{3}}\bigg)|\xi|,\eeq where $k\in\{1,2,3,4\}$.

Next, we estimate
$$I_{11}=\int_{\mathbb{B}^{4}}\frac{(1-|y|^{2})}{[x,y]^{2}}dV(y),~I_{12}=\int_{\mathbb{B}^{4}}\frac{(x\otimes
y)}{|x-y|^{2}[x,y]^{2}}dV(y)$$ and
$$I_{13}=\int_{\mathbb{B}^{4}}\left(\frac{(1-|y|^{2})}{[x,y]^{2}}+\frac{x\otimes
y}{[x,y]^{3}}\right)dV(y).$$

First, we know from elementary calculations  that

\beq\label{0-n=4}
I_{11}&\leq&\int_{\mathbb{B}^{4}}\frac{(1-|x|^{2}|y|^{2})}{[x,y]^{2}}dV(y)=
A_{3}\int_{0}^{1}\rho^{3}\left(\int_{\mathbb{S}^{3}}\frac{1-|x|^{2}\rho^{2}}{|x\rho-\zeta|^{2}}d\sigma(\zeta)\right)d\rho\\
\nonumber&=&\frac{A_{3}}{4}. \eeq

Let $z=\phi_{x}(y)\in\Aut(\mathbb{B}^{4})$. Then, by (\ref{II}),
(\ref{III}), (\ref{eq-j-1}) and (\ref{c-20v}), we obtain

\beq\label{1-n=4} I_{12}&=&(1-|x|^{2})^{2}\int_{\mathbb{B}^{4}}\frac{(1-|z|^{2})}{|z|^{2}[x,z]^{6}}dV(z)\\
\nonumber
&=&A_{3}(1-|x|^{2})^{2}\int_{0}^{1}\left(\int_{\mathbb{S}^{3}}\frac{\rho(1-\rho^{2})}{|x\rho-\zeta|^{6}}d\sigma(\zeta)\right)d\rho\\
\nonumber
&=&A_{3}(1-|x|^{2})^{2}\int_{0}^{1}\frac{\rho(1-\rho^{2})}{(1-|x|^{2}\rho^{2})^{3}}d\rho\\
\nonumber
&\leq&A_{3}(1-|x|^{2})^{2}\int_{0}^{1}\frac{\rho}{(1-|x|^{2}\rho^{2})^{2}}d\rho\leq\frac{A_{3}}{2}.
\eeq

Since \beqq \frac{1-|y|^{2}}{[x,y]^{2}}+\frac{x\otimes
y}{[x,y]^{3}}&=&(1-|y|^{2})\left(\frac{1}{[x,y]^{2}}+\frac{1-|x|^{2}}{[x,y]^{2}(1-|x|)}\right)\\
&\leq&3\frac{(1-|x|^{2}|y|^{2})}{[x,y]^{2}}, \eeqq we see that

\beq\label{2-n=4}I_{13}&\leq&3A_{3}\int_{0}^{1}\rho^{3}\left(\int_{\mathbb{S}^{3}}\frac{1-|x|^{2}\rho^{2}}{|x\rho-\zeta|^{2}}d\sigma(\zeta)\right)d\rho
 =\frac{3}{4}A_{3}. \eeq

It follows from  (\ref{eq-hj-s1}), (\ref{0-n=4}), (\ref{1-n=4}),
(\ref{2-n=4}) and Lemma \Ref{K-3} that

\beq\label{eq-hj-2g}|D_{G[g]}(x)\xi|&\leq&\|g\|_{\infty}\int_{\mathbb{B}^{4}}\left|\langle\nabla
G_{2,4}(x,y),\xi\rangle\right| dV(y)\\ \nonumber
&\leq&2|c_{4}|\big(I_{11}+2I_{12}+I_{13}\big)|\xi|\leq2|c_{4}|\left(\frac{A_{3}}{4}+A_{3}+\frac{3A_{3}}{4}\right)|\xi|\\
\nonumber &=&\frac{1}{2}|\xi|.\eeq

Therefore, by (\ref{eq-ca-g1}),  (\ref{eq-hj-87}) and
(\ref{eq-hj-2g}), we conclude that there exists a positive constant
$M_{5}$ such that

\be\label{eq-ca-h24}|D_{G[g]}(x)\xi|\leq M_{5}\|g\|_{\infty}\xi.\ee

\bst\label{bst-4.0} The Lipschitz continuity of $f$. \est

By (\ref{eq-ca-k8}), (\ref{eq-c-h10}) and  (\ref{eq-ca-h24}), we see
that there is a  constant
$M_{6}(n,\|\varphi_{2}\|_{\infty},\|g\|_{\infty})$ such that
$$|D_{f}(x)|\leq M_{6}(n,\|\varphi_{2}\|_{\infty},\|g\|_{\infty}),$$ which yields that, for any
$a,b\in\mathbb{B}^{n}$,
\beqq|f(a)-f(b)|&=&\left|\int_{[a,b]}D_{f}(x)dx\right|=\left|\int_{[a,b]}D_{f}(x)\frac{dx}{|dx|}|dx|\right|\leq\int_{[a,b]}|D_{f}(x)||dx|\\
&\leq& M_{6}(n,\|\varphi_{2}\|_{\infty},\|g\|_{\infty})|a-b|\\
&=&M_{6}(n,\|\varphi_{2}\|_{\infty},\|g\|_{\infty})\left(\frac{|a-b|}{\omega(|a-b|)}\right)\omega(|a-b|)\\
&\leq&M_{6}(n,\|\varphi_{2}\|_{\infty},\|g\|_{\infty})\frac{2}{\omega(2)}\omega(|a-b|).\eeqq
The proof of this theorem is complete. \qed

\bigskip

{\bf Acknowledgements:}    This research was partly supported by the
Hunan Provincial Education Department Outstanding Youth Project (No.
18B365), the Science and Technology Plan Project of Hengyang City
(No. 2018KJ125),  the Science and Technology Plan Project of Hunan
Province (No. 2016TP1020), the Science and Technology Plan Project
of Hengyang City (No. 2017KJ183), and the Application-Oriented
Characterized Disciplines, Double First-Class University Project of
Hunan Province (Xiangjiaotong [2018]469).


\end{document}